# TWO-PHASE SAMPLING IN ESTIMATION OF POPULATION MEAN IN THE PRESENCE OF NON-RESPONSE

Manoj Kr. Chaudhary, Anil Prajapati and †Rajesh Singh

Department of StatisticsBanaras Hindu University

Varanasi-221005

†Corresponding author

**Abstract**

The present paper presents the detail discussion on estimation of population mean in simple random sampling in the presence of non-response. Motivated by Gupta and Shabbir (2008), we have suggested the class of estimators of population mean using an auxiliary variable under non-response. A theoretical study is carried out using two-phase sampling scheme when the population mean of auxiliary variable is not known. An empirical study has also been done in the support of theoretical results.

## 1. Introduction

The auxiliary information is generally used to improve the efficiency of the estimators. Cochran (1940) proposed the ratio estimator for estimating the population mean whenever study variable is positively correlated with auxiliary variable. Contrary to the situation of ratio estimator, if the study and auxiliary variables are negatively correlated, Murthy (1964) suggested the product estimator to estimate the population mean. Hansen et al. (1953) proposed the difference estimator which was subsequently modified to provide the linear regression estimator for the population mean or total. Mohanty (1967) suggested an estimator by combining the ratio and regression methods for estimating the population parameters. In order to estimate the population mean or population total of the study character utilizing auxiliary information, several other authors including Srivastava ( 1971), Reddy (1974), Ray and Sahai (1980), Srivenkataramana (1980), Srivastava and Jhajj (1981) and Singh and Kumar (2008, 2011) have proposed estimators which lead improvements over usual per unit estimator.

It is observed that the non-response is a common problem in any type of survey. Hansen and Hurwitz (1946) were the first to contract the problem of non-response while conducting mail surveys. They suggested a technique, known as 'sub-sampling of non-respondents', to deal with the problem of non-response and its adjustments. In fact they developed an unbiased estimator for population mean in the presence of non-response by dividing the population into two groups, viz. response group and non-response group. To avoid bias due to non-response, they suggested for taking a sub-sample of the non-responding units.

Let us consider a population consists of $N$ units and a sample of size $n$ is selected from the population using simple random sampling without replacement (SRSWOR) scheme. Let us assume that $Y$ and $X$ be the study and auxiliary variables with respective population means $\overline{Y}$ and $\overline{X}$. Let us consider the situation in which study variable is subjected to non-response and auxiliary variable is free from the non-response. It is observed that there are $n_1$ respondent and $n_2$ non-respondent units in the sample of $n$ units for the study variable. Using the technique of sub sampling of non-respondents suggested by Hansen and Hurwitz (1946), we select a sub-sample of $h_2$ non-respondent units from $n_2$ units such that $h_2 = n_2/k, k \geq 1$ and collect the information on sub-sample by personal interview method. The usual sample mean, ratio and regression estimators for estimating the population mean $\overline{Y}$ under non-response are respectively represented by

$$\overline{y}^* = \frac{n_1 \overline{y}_{n1} + n_2 \overline{y}_{h2}}{n} \tag{1.1}$$

$$\bar{y}_R^* = \frac{\bar{y}^*}{\bar{x}}\bar{X} \tag{1.2}$$

$$\bar{y}_{lr}^* = \bar{y}^* + b\left(\bar{X} - \bar{x}\right) \tag{1.3}$$

where $\bar{y}_{n1}$ and $\bar{y}_{h2}$ are the means based on $n_1$ respondent and $h_2$ non-respondent units respectively. $\bar{x}$ is the sample mean estimator of population mean $\bar{X}$, based on sample of size $n$ and $b$ is the sample regression coefficient of $Y$ on $X$.

The variance and mean square errors (MSE) of the above estimators $\bar{y}^*$, $\bar{y}_R^*$ and $\bar{y}_{lr}^*$ are respectively given by

$$V\left(\bar{y}^*\right) = \left(\frac{1}{n} - \frac{1}{N}\right)S_Y^2 + \frac{(k-1)}{n}W_2 S_{Y2}^2 \tag{1.4}$$

$$MSE\left(\bar{y}_R^*\right) = \left(\frac{1}{n} - \frac{1}{N}\right)\bar{Y}^2\left(C_Y^2 + C_X^2 - 2\rho C_X C_Y\right) + \frac{(k-1)}{n}W_2 S_{Y2}^2 \tag{1.5}$$

$$MSE\left(\bar{y}_{lr}^*\right) = \left(\frac{1}{n} - \frac{1}{N}\right)\bar{Y}^2 C_Y^2\left(1 - \rho^2\right) + \frac{(k-1)}{n}W_2 S_{Y2}^2 \tag{1.6}$$

where $S_Y^2$ and $S_X^2$ are respectively the mean squares of $Y$ and $X$ in the population. $C_Y\left(= S_Y/\bar{Y}\right)$ and $C_X\left(= S_X/\bar{X}\right)$ are the coefficients of variation of $Y$ and $X$ respectively. $S_{Y2}^2$ and $W_2$ are respectively the mean square and non-response rate of the non-response group in the population for the study variable $Y$. $\rho$ is the population correlation coefficient between $Y$ and $X$.

When the information on population mean of auxiliary variable is not available, one can use the two-phase sampling scheme in obtaining the improved estimator rather than the previous ones. Neyman (1938) was the first who gave concept of two-phase sampling in estimating the population parameters. Two-phase sampling is cost effective as well as easier. This sampling scheme is used to obtain the information about auxiliary variable cheaply from a bigger sample at first phase and relatively small sample at the second stage. Sukhatme (1962) used two-phase sampling scheme to propose a general ratio-type estimator. Rao (1973) used two-phase sampling to stratification, non-response problems and investigative comparisons. Cochran (1977) supplied some basic information for two-phase sampling. Sahoo et al. (1993) provided regression approach in estimation by using two auxiliary variables for two-phase sampling. In the sequence of improving the efficiency of the estimators, Singh and Upadhyaya (1995) suggested a generalized estimator to estimate population mean using two auxiliary variables in two-phase sampling.

In estimating the population mean $\bar{Y}$, if $\bar{X}$ is unknown, first, we obtain the estimate of it using two-phase sampling scheme and then estimate $\bar{Y}$. Under two-phase sampling scheme, first we select a larger sample of $n'$ units from the population of size $N$ with the help of SRSWOR scheme. Secondly, we select a small sample of size $n$ from $n'$ units. Let us again assume that the situation in which the non-response is observed on study variable only and auxiliary variable is free from the non-response. The usual ratio and regression estimators of population mean $\bar{Y}$ under two-phase sampling in the presence of non-response are respectively given by

$$\bar{y}_R^{**} = \frac{\bar{y}^*}{\bar{x}} \bar{x}' \tag{1.7}$$

and $$\bar{y}_{lr}^{**} = \bar{y}^* + b\left(\bar{x}' - \bar{x}\right) \tag{1.8}$$

where $\bar{x}'$ is the mean based on $n'$ units for the auxiliary variable.

The MSE's of the estimators $\bar{y}_R^{**}$ and $\bar{y}_{lr}^{**}$ are respectively represented by the following expressions

$$MSE\left(\bar{y}_R^{**}\right) = \bar{Y}^2\left[\left(\frac{1}{n'} - \frac{1}{N}\right)C_Y^2 + \left(\frac{1}{n} - \frac{1}{n'}\right)\left(C_Y^2 + C_X^2 - 2\rho C_X C_Y\right)\right] + \frac{(k-1)}{n} W_2 S_{Y2}^2 \tag{1.9}$$

and

$$MSE\left(\bar{y}_{lr}^{**}\right)=\bar{Y}^2\left[\left(\frac{1}{n'}-\frac{1}{N}\right)C_Y^2+\left(\frac{1}{n}-\frac{1}{n'}\right)C_Y^2\left(1-\rho^2\right)\right]+\frac{(k-1)}{n}W_2S_{Y2}^2 \tag{1.10}$$

In the present paper, we have discussed the study of non-response of a general class of estimators using an auxiliary variable. We have suggested the class of estimators in two-phase sampling when the population mean of auxiliary variable is unknown. The optimum property of the class is also discussed and it is compared to ratio and regression estimators under non-response. The theoretical study is also supported with the numerical illustrations.

## 2. Suggested Class of Estimators

Let us assume that the non-response is observed on the study variable and auxiliary variable provides complete response on the units. Motivated by Gupta and Shabbir (2008), we suggest a class of estimators of population mean $\bar{Y}$ under non-response as

$$\bar{y}_t^* = \left[\alpha_1\bar{y}^* + \alpha_2\left(\bar{X}-\bar{x}\right)\right]\left(\frac{\eta\bar{X}+\lambda}{\eta\bar{x}+\lambda}\right) \tag{2.1}$$

where $\alpha_1$ and $\alpha_2$ are the constants and whose values are to be determined. $\lambda$ and $\eta(\neq 0)$ are either constants or functions of the known parameters.

In order to obtain the bias and MSE of $\bar{y}_t^*$, we use the large sample approximation. Let us assume that

$\bar{y}^* = \bar{Y}(1+e_1)$, $\bar{x} = \bar{X}(1+e_2)$

such that $E(e_1)=E(e_2)=0$,

$$E\left(e_1^2\right)=\frac{V\left(\bar{y}^*\right)}{\bar{Y}^2}=\left(\frac{1}{n}-\frac{1}{N}\right)C_Y^2+\frac{(k-1)}{n}W_2\frac{S_{Y2}^2}{\bar{Y}^2},$$

$$E\left(e_2^2\right)=\frac{V\left(\bar{x}\right)}{\bar{X}^2}=\left(\frac{1}{n}-\frac{1}{N}\right)C_X^2$$

and $E(e_1e_2)=\dfrac{Cov\left(\bar{y}^*,\bar{x}\right)}{\bar{Y}\bar{X}}=\left(\dfrac{1}{n}-\dfrac{1}{N}\right)\rho C_X C_Y$.

Putting the values of $\bar{y}^*$ and $\bar{x}$ form the above assumptions in the equation (2.1), we get

$$\bar{y}_t^* - \bar{Y} \cong \bar{Y}(\alpha_1 - 1) + \alpha_1 \bar{Y}\left(e_1 - \tau e_2 - \tau e_1 e_2 + \tau^2 e_2^2\right) - \alpha_2 \bar{X}\left(e_2 - \tau e_2^2\right) \tag{2.2}$$

On taking expectation of the equation (2.2), the bias of $\bar{y}_t^*$ to the first order of approximation is given by

$$B\left(\bar{y}_t^*\right) = E\left(\bar{y}_t^* - \bar{Y}\right) = \bar{Y}(\alpha_1 - 1) + \left(\frac{1}{n} - \frac{1}{N}\right)\left[\alpha_1 \bar{Y}\left(\tau^2 C_X^2 - \tau\rho C_X C_Y\right) + \alpha_2 \bar{X}\tau C_X^2\right] \tag{2.3}$$

Squaring both the sides of the equation (2.2) and taking expectation, we can obtain the MSE of $\bar{y}_t^*$ to the first order of approximation as

$$MSE\left(\bar{y}_t^*\right) = \bar{Y}^2(\alpha_1 - 1)^2 + \left(\frac{1}{n} - \frac{1}{N}\right)\Big[\alpha_1^2 \bar{Y}^2\left(C_Y^2 + \tau^2 C_X^2 - 2\tau\rho C_X C_Y\right) + \alpha_2^2 \bar{X}^2 C_X^2$$
$$- 2\alpha_1\alpha_2 \bar{Y}\bar{X} C_X\left(\rho C_Y - \tau C_X\right)\Big] + \alpha_1^2 \frac{(k-1)}{n} W_2 S_{Y2}^2 \tag{2.4}$$

In the sequence of obtaining the best estimator within the suggested class with respect to $\alpha_1$ and $\alpha_2$, we obtain the optimum values of $\alpha_1$ and $\alpha_2$. On differentiating $MSE\left(\bar{y}_t^*\right)$ with respect to $\alpha_1$ and $\alpha_2$ and equating the derivatives to zero, we have

$$\frac{\partial MSE\left(\bar{y}_t^*\right)}{\partial \alpha_1} = \bar{Y}^2(\alpha_1 - 1) + \left(\frac{1}{n} - \frac{1}{N}\right)\left[\bar{Y}^2 \alpha_1\left(C_Y^2 + \tau^2 C_X^2 - 2\tau\rho C_X C_Y\right) - \alpha_2 \bar{X}\bar{Y} C_X\left(\rho C_Y - \tau C_X\right)\right]$$
$$+ \alpha_1 \frac{(k-1)}{n} W_2 S_{Y2}^2 = 0 \tag{2.5}$$

$$\frac{\partial MSE\left(\bar{y}_t^*\right)}{\partial \alpha_2} = \left(\frac{1}{n} - \frac{1}{N}\right)\left[\alpha_2 \bar{X}^2 C_X^2 - \alpha_1 \bar{X}\bar{Y} C_X\left(\rho C_Y - \tau C_X\right)\right] = 0 \tag{2.6}$$

Solving the equations (2.4) and (2.5), we get

$$\alpha_1(opt) = \frac{1}{1 + \left(\frac{1}{n} - \frac{1}{N}\right) C_Y^2\left(1 - \rho^2\right) + \frac{(k-1)}{n} W_2 \frac{S_{Y2}^2}{\bar{Y}^2}} \tag{2.7}$$

and $\alpha_2(\text{opt}) = \dfrac{\alpha_1(\text{opt})\overline{Y}(\rho C_Y - \tau C_X)}{\overline{X}C_X}$ (2.8)

Substituting the values of $\alpha_1(\text{opt})$ and $\alpha_2(\text{opt})$ from equations (2.7) and (2.8) into the equation (2.4), the MSE of $\overline{y}_t^*$ is given by the following expression.

$$\text{MSE}\left(\overline{y}_t^*\right)_{\min} = \frac{\text{MSE}\left(\overline{y}_{lr}^*\right)}{1 + \left(\frac{1}{n} - \frac{1}{N}\right)C_Y^2\left(1-\rho^2\right) + \frac{(k-1)}{n}W_2\frac{S_{Y2}^2}{\overline{Y}^2}} \tag{2.9}$$

## 3. Suggested Class in Two-Phase Sampling

It is generally seen that the population mean of auxiliary variable, $\overline{X}$ is not known. In this situation, we may use the two-phase sampling scheme to find out the estimate of $\overline{X}$. Using two-phase sampling, we now suggest a class of estimators of population mean $\overline{Y}$ in the presence of non-response when $\overline{X}$ is unknown, as

$$\overline{y}_t^{**} = \left[\alpha_1\overline{y}^* + \alpha_2\left(\overline{x}' - \overline{x}\right)\right]\left(\frac{\eta\overline{x}' + \lambda}{\eta\overline{x} + \lambda}\right) \tag{3.1}$$

### 3.1 Bias and MSE of $\overline{y}_t^{**}$

By applying the large sample approximation, we can obtain the bias and mean square error of $\overline{y}_t^{**}$. Let us assume that

$$\overline{y}^* = \overline{Y}(1+e_1),\ \overline{x} = \overline{X}(1+e_2) \text{ and } \overline{x}' = \overline{X}(1+e_3)$$

such that $E(e_1) = E(e_2) = E(e_3) = 0$,

$E(e_1^2) = \left(\frac{1}{n} - \frac{1}{N}\right)C_Y^2 + \frac{(k-1)}{n}W_2\frac{S_{Y2}^2}{\overline{Y}^2}$ , $E(e_2^2) = \left(\frac{1}{n} - \frac{1}{N}\right)C_X^2$,

$E(e_3^2) = \left(\frac{1}{n'} - \frac{1}{N}\right)C_X^2$ , $E(e_1e_2) = \left(\frac{1}{n} - \frac{1}{N}\right)\rho C_X C_Y$,

$$E(e_1 e_3) = \left(\frac{1}{n'} - \frac{1}{N}\right)\rho C_X C_Y \text{ and } E(e_2 e_3) = \left(\frac{1}{n'} - \frac{1}{N}\right)C_X^2 .$$

Under the above assumption, the equation (3.1) gives

$$\bar{y}_t^{**} - \bar{Y} = \bar{Y}(\alpha_1 - 1) + \alpha_1 \bar{Y}\left(e_1 + \tau^2 e_2^2 - \tau e_2 - \tau e_1 e_2 + \tau e_3 + \tau e_1 e_3 - \tau^2 e_2 e_3\right)$$

$$+ \alpha_2 \bar{X}\left(e_3 - e_2 \tau e_2 e_3 + \tau e_2^2 + \tau e_3^2 - \tau e_2 e_3\right) \tag{3.2}$$

Taking expectation of both the sides of equation (3.2), we get the bias of $\bar{y}_t^{**}$ up to the first order of approximation as

$$B\left(\bar{y}_t^{**}\right) = \bar{Y}(\alpha_1 - 1) + \left(\frac{1}{n} - \frac{1}{n'}\right)\tau\left[\alpha_1 \bar{Y}\left(C_X^2 - \rho C_X C_Y\right) + \alpha_2 \bar{X} C_X^2\right] \tag{3.3}$$

The MSE of $\bar{y}_t^{**}$ up to the first order of approximation can be obtained by the following expression

$$MSE\left(\bar{y}_t^{**}\right) = E\left(\bar{y}_t^{**} - \bar{Y}\right)^2 = \bar{Y}^2 (\alpha_1 - 1)^2$$

$$+ \alpha_1^2 \bar{Y}^2 \left[\left(\frac{1}{n} - \frac{1}{N}\right)C_Y^2 + \left(\frac{1}{n} - \frac{1}{n'}\right)\left(\tau^2 C_X^2 - 2\tau\rho C_X C_Y\right) + \frac{(k-1)}{n} W_2 \frac{S_{Y2}^2}{\bar{Y}^2}\right]$$

$$+ \left(\frac{1}{n} - \frac{1}{n'}\right)\left[\alpha_2^2 \bar{X}^2 C_X^2 + 2\alpha_1 \alpha_2 \bar{X}\bar{Y}\left(\tau C_X^2 - \rho C_X C_Y\right)\right] \tag{3.4}$$

**3.2 Optimum Values of $\alpha_1$ and $\alpha_2$**

On differentiating $MSE\left(\bar{y}_t^{**}\right)$ with respect to $\alpha_1$ and $\alpha_2$ and equating the derivatives to zero, we get the normal equations

$$\frac{\partial MSE\left(\bar{y}_t^{**}\right)}{\partial \alpha_1} = \bar{Y}^2 (\alpha_1 - 1) + \alpha_1 \bar{Y}^2 \left[\left(\frac{1}{n} - \frac{1}{N}\right)C_Y^2 + \left(\frac{1}{n} - \frac{1}{n'}\right)\left(\tau^2 C_X^2 - 2\tau\rho C_X C_Y\right) + \frac{(k-1)}{n} W_2 \frac{S_{Y2}^2}{\bar{Y}^2}\right]$$

$$+ \left(\frac{1}{n} - \frac{1}{n'}\right)\alpha_2 \bar{X}\bar{Y}\left(\tau C_X^2 - \rho C_X C_Y\right) = 0 \tag{3.5}$$

and $\dfrac{\partial MSE\left(\bar{y}_t^{**}\right)}{\partial \alpha_2} = \left(\dfrac{1}{n} - \dfrac{1}{n'}\right)\left[\alpha_2 \bar{X}^2 C_X^2 + \alpha_1 \bar{X}\bar{Y}\left(\tau C_X^2 - \rho C_X C_Y\right)\right] = 0$ (3.6)

From equations (3.5) and (3.6), we get the optimum values of $\alpha_1$ and $\alpha_2$ as

$$\alpha_1(opt) = \frac{1}{1 + \left(\frac{1}{n} - \frac{1}{N}\right)C_Y^2 - \left(\frac{1}{n} - \frac{1}{n'}\right)\rho^2 C_Y^2 + \frac{(k-1)}{n} W_2 \frac{S_{Y2}^2}{\bar{Y}^2}} \tag{3.7}$$

and $\alpha_2(opt) = \dfrac{\alpha_1(opt)\bar{Y}\left(\rho C_Y - \tau C_X\right)}{\bar{X} C_X}$ (3.8)

On substituting the optimum values of $\alpha_1$ and $\alpha_2$, the equation (3.4) provides

minimum MSE of $\bar{y}_t^{**}$

$$MSE\left(\bar{y}_t^{**}\right)_{min} = \frac{MSE\left(\bar{y}_{lr}^{**}\right)}{1 + \left(\frac{1}{n} - \frac{1}{N}\right)C_Y^2 - \left(\frac{1}{n} - \frac{1}{n'}\right)\rho^2 C_Y^2 + \frac{(k-1)}{n} W_2 \frac{S_{Y2}^2}{\bar{Y}^2}} \tag{3.9}$$

## 4. Empirical Study

In the support of theoretical results, some numerical illustrations are given below:

**4.1** In this section, we have illustrated the relative efficiency of the estimators $\bar{y}_R^*$, $\bar{y}_{lr}^*$ and $\bar{y}_t^*(opt)$ with respect to $\bar{y}^*$. For this purpose, we have considered the data used by Kadilar and Cingi (2006). The details of the population are given below:

$N = 200$, $n = 50$, $\bar{Y} = 500$, $\bar{X} = 25$, $C_Y = 15$, $C_X = 2$, $\rho = 0.90$

$k = 1.5$, $S_{Y2}^2 = \frac{4}{5} S_Y^2$

**Table 1. Percentage Relative Efficiency (PRE) with respect to $\bar{y}^*$**

| $W_2$ | Estimator | | |
|---|---|---|---|
| | $\bar{y}_R^*$ | $\bar{y}_{lr}^*$ | $\bar{y}_t^*(opt)$ |
| 0.1 | 126.74 | 432.88 | 788.38 |
| 0.2 | 125.13 | 373.03 | 746.53 |
| 0.3 | 123.70 | 331.43 | 722.93 |
| 0.4 | 122.42 | 300.83 | 710.33 |
| 0.5 | 121.28 | 277.37 | 704.87 |

**4.2** The present section presents the relative efficiency of the estimators $\bar{y}_R^{**}$, $\bar{y}_{lr}^{**}$ and $\bar{y}_t^{**}(opt)$ with respect to $\bar{y}^*$. There are two data sets which have been considered to illustrate the theoretical results.

**Data Set 1:**

The population considered by Srivastava (1993) is used to give the numerical interpretation of the present study. The population of seventy villages in a Tehsil of India along with their cultivated area (in acres) in 1981 is considered. The cultivated area (in acres) is taken as study variable and the population is assumed to be auxiliary variable. The population parameters are given below:

$N = 70$, $n' = 40$, $n = 25$, $\bar{Y} = 981.29$, $\bar{X} = 1755.53$, $S_Y = 613.66$,

$S_X = 1406.13$, $S_{Y2} = 244.11$, $\rho = 0.778$, $k = 1.5$

**Table 2. Percentage Relative Efficiency with respect to $\bar{y}^*$**

| $W_2$ | Estimator | | |
|---|---|---|---|
| | $\bar{y}_R^*$ | $\bar{y}_{lr}^*$ | $\bar{y}_t^*(opt)$ |
| 0.1 | 125.48657 | 153.56020 | 154.57983 |
| 0.2 | 125.10358 | 152.57858 | 153.60848 |
| 0.3 | 124.73193 | 151.63228 | 152.67552 |
| 0.4 | 124.37111 | 150.71945 | 151.77449 |
| 0.5 | 124.02068 | 149.83834 | 150.90579 |

**Data Set 2:**

Now, we have used another population considered by Khare and Sinha (2004). The data are based on the physical growth of upper-socio-economic group of 95 school children of Varanasi district under an ICMR study, Department of Paediatrics, Banaras Hindu University, India during 1983-84. The details are given below:

$N = 95,\ n' = 70,\ \ n = 35,\ \bar{Y} = 19.4968,\ \bar{X} = 55.8611,\ S_Y = 3.0435,\ S_X = 3.2735,$

$S_{Y2} = 2.3552,\ \rho = 0.8460,\ k = 1.5$

**Table 3. Percentage Relative Efficiency with respect to $\bar{y}^*$**

| $W_2$ | Estimator | | |
|---|---|---|---|
| | $\bar{y}_R^*$ | $\bar{y}_{lr}^*$ | $\bar{y}_t^*(opt)$ |
| 0.1 | 159.61889 | 217.83004 | 217.99278 |
| 0.2 | 155.61224 | 207.27149 | 207.43596 |
| 0.3 | 152.10325 | 198.44091 | 198.58540 |
| 0.4 | 149.01829 | 190.94488 | 190.94488 |
| 0.5 | 146.26158 | 184.51722 | 184.66554 |

## 5. Conclusion

The study of a general class of estimators of population mean under non-response has been presented. We have also suggested a class of estimators of population mean in the presence of non-response using two-phase sampling when population mean of auxiliary variable is not known. The optimum property of the suggested class has been discussed. We have compared the optimum estimator with some existing estimators through numerical study. The Tables 1, 2 and 3 represent the percentage relative efficiency of the optimum estimator of suggested class, linear regression estimator and ratio estimator with respect to sample mean estimator. In the above tables, we have observed that the percentage relative efficiency of the optimum estimator is higher than the linear regression and ratio estimators. It is also observed that the percentage relative efficiency decreases with increase in non-response.